\documentclass[aos,MSNbibl,nameyear,dvips]{arximspdf}
\usepackage{graphicx}

\doi{10.1214/12-AOS994} 
\volume{40}
\issue{2}
\pubyear{2012}
\firstpage{941}
\lastpage{963}

\makeatletter
\let\hat\widehat
\let\tilde\widetilde
\newcommand{\argmax}{\operatorname{\arg\max}}

\newcommand\cA{{\mathcal A}}
\newcommand\cG{{\mathcal G}}
\newcommand\cK{{\mathcal K}}
\newcommand\cM{{\mathcal M}}
\newcommand\cQ{{\mathcal Q}}
\newcommand\intersect{\cap}
\newcommand{\norm}[1]{\|#1\|}
\newcommand{\Set}[1]{\{#1\}}
\newcommand\R{\mathbb{R}}
\newcommand\E{\mathbb{E}}
\renewcommand\P{\mathbb{P}}

\newtheorem{theorem}{Theorem}
\newtheorem{lemma}[theorem]{Lemma}
\newproclaim{remark}{Remark}
\makeatother

\begin{document}
\begin{frontmatter}

\title{Manifold estimation and singular deconvolution under Hausdorff loss}
\runtitle{Manifold estimation}

\begin{aug}
\author[A]{\fnms{Christopher R.} \snm{Genovese}\thanksref{aut1}\ead[label=e1]{genovese@stat.cmu.edu}},
\author[B]{\fnms{Marco} \snm{Perone-Pacifico}\thanksref{aut2}\ead[label=e2]{marco.peronepacifico@uniroma1.it}},\\
\author[C]{\fnms{Isabella} \snm{Verdinelli}\thanksref{aut2}\ead[label=e3]{isabella@stat.cmu.edu}}
\and
\author[A]{\fnms{Larry} \snm{Wasserman}\thanksref{aut3}\corref{}\ead[label=e4]{larry@stat.cmu.edu}}
\thankstext{aut1}{Supported by NSF Grant DMS-08-06009.}
\thankstext{aut2}{Supported by Italian National Research Grant PRIN
2008.}
\thankstext{aut3}{Supported by NSF Grant DMS-08-06009, Air Force Grant
FA95500910373
and Sapienza University of Rome grant for visiting professors 2009.}
\runauthor{Genovese, Perone-Pacifico, Verdinelli and Wasserman}
\affiliation{Carnegie Mellon University, Sapienza University of Rome,
Carnegie Mellon University and Sapienza University of Rome, and
Carnegie Mellon University}
\address[A]{C. R. Genovese\\L. Wasserman\\Department of Statistics\\Carnegie Mellon
University\\
Pittsburgh, Pennsylvania 15213\\USA\\ \printead{e1}\\\phantom{E-mail: }\printead*{e4}} 
\address[B]{M. Perone-Pacifico\\Department of Statistical Sciences\\
Sapienza University of Rome\\
Rome\\ Italy\\ \printead{e2}}
\address[C]{I. Verdinelli\\Department of Statistics\\Carnegie Mellon
University\\
Pittsburgh, Pennsylvania 15213\\USA\\and\\
Department of Statistical Sciences\\
Sapienza University of Rome\\Rome\\ Italy\\ \printead{e3}}
\end{aug}

\received{\smonth{9} \syear{2011}}
\revised{\smonth{1} \syear{2012}}

%
\begin{abstract}
We find lower and upper bounds
for the risk of estimating a~manifold
in Hausdorff distance
under
several models.
We also show that there are close connections between manifold estimation
and the problem of deconvolving a~singular measure.
\end{abstract}

%

\end{frontmatter}

\section{Introduction}

Manifold learning is an area of intense research activity in machine
learning and statistics. Yet a~very basic question about manifold
learning is still
open, namely, how well can we estimate a~manifold from~$n$ noisy
samples? In this paper we investigate this question under various
assumptions.

Suppose we observe a~random sample $Y_1,\ldots, Y_n \in\mathbb{R}^D$
that lies on or near
a~$d$-manifold $M$ where $d< D$.
The question we address is:
what is the minimax risk under Hausdorff distance for estimating $M$?
Our main assumption is that
$M$ is a~$d$-dimensional, smooth Riemannian submanifold in $\mathbb{R}^D$;
the precise conditions on $M$ are given in
Section \ref{sec::conditions}.

Let $Q$ denote the distribution of $Y_i$.
We shall see that $Q$ depends on several things,
including the manifold $M$, a~distribution $G$ supported on $M$ and
a~model for the noise.
We consider three noise models.
The first is the \textit{noiseless} model in which
$Y_1,\ldots, Y_n$ is a~random sample from $G$.
The second is the \textit{clutter noise} model, in which
%
\begin{equation}
Y_1,\ldots, Y_n \sim(1-\pi)U + \pi G,
\end{equation}
where $U$ is a~uniform distribution on
a~compact set $\cK\subset\mathbb{R}^D$
with nonempty interior,
and
$G$ is supported on $M$.
(When $\pi=1$ we recover the noiseless case.)
The third is the \textit{additive} model,
%
\begin{equation}
Y_i = X_i + Z_i,
\end{equation}
where
$X_1,\ldots, X_n \sim G$,
$G$ is supported on $M$,
and
the noise variables
$Z_1, \ldots,\allowbreak Z_n$
are a~sample from a~distribution
$\Phi$ on $\mathbb{R}^D$ which we take to be Gaussian.
In this case, the distribution $Q$ of $Y$
is a~convolution of $G$ and $\Phi$ written
$Q = G\star\Phi$.

In a~previous paper
[\citet{us::2010}],
we considered a~noise model
in which the noise is perpendicular to the manifold.
This model is also considered in \citet{smale2}.
Since we have already studied that model,
we shall not consider it further here.

In the additive model,
estimating $M$ is related to estimating the distribution~$G$,
a~problem that is usually called \textit{deconvolution}
[\citet{Fan}].
The problem of deconvolution is well studied in the statistical literature,
but in the manifold case there is an interesting complication: the
measure $G$
is singular because it puts all its mass on a~subset of $\mathbb{R}^D$
that has zero Lebesgue measure
(since the manifold has dimension $d<D$).
Deconvolution of singular measures has not received
as much attention as standard deconvolution problems
and raises interesting challenges.

Each noise model gives rise to
a~class of distributions
$\cQ$ for $Y$
defined more precisely in Section
\ref{sec::conditions}.
We are interested in the
minimax risk
%
\begin{equation}
R_n \equiv R_n(\cQ) = \inf_{\hat M}\sup_{Q\in\cQ}\E_Q [ H(\hat M,M)],
\end{equation}
where the infimum is over all estimators $\hat M$,
and $H$ is the Hausdorff distance
[defined in equation (\ref{eq::hausdorff})].
Note that finding the minimax risk is equivalent to finding
the \textit{sample complexity}
$n(\varepsilon) = \inf \{n\dvtx R_n \leq\varepsilon \}.$
We emphasize that the goal of this paper is to find the minimax rates,
not to find practical estimators.
We use the Hausdorff distance because
it is one of the most commonly used metrics for
assessing the accuracy of set-valued estimators.
One could of course create other loss functions
and study their properties, but this is beyond the scope of this paper.
Finally, we remark that our upper bounds
sometimes differ from our lower bounds
by a~logarithmic factor.
This is a~common phenomenon when dealing
with Hausdorff distance (and sup norm in function estimation problems).
Currently, we do not know how to eliminate the log factor.

\subsection{Related work}
In the additive noise case,
estimating a~manifold is related to deconvolution problems
such as those in
\citet{Fan}, \citet{FanTruong}
and
\citet{Stefanski1990229}.
More closely related is the problem of estimating the support of a~distribution
in the presence of noise as discussed, for example, in
\citet{Meister20061702}.\vadjust{\goodbreak}

There is a~vast literature on manifold estimation.
Much of the literature deals with
using manifolds for the purpose of dimension reduction.
See, for example,
\citet{baraniuk} and references therein.
We are interested instead in actually estimating the manifold itself.
There is a~literature on this problem in
the field of computational geometry;
see \citet{Dey}.
However, very few papers allow for noise in the statistical sense,
by which we mean observations drawn randomly
from a~distribution.
In the literature on computational geometry,
observations are called noisy if they depart from the underlying
manifold in a~very specific way:
the observations have to be close to the manifold but not too close
to each other.
This notion of noise is quite different from random sampling
from a~distribution.
An exception is \citet{smale},
who constructed the following estimator:
Let
$I = \{ i\dvtx \hat{p}(Y_i) > \lambda\}$
where
$\hat{p}$ is a~density estimator.
They define
$\hat M = \bigcup_{i\in I} B_D(Y_i,\varepsilon)$
where $B_D(Y_i,\varepsilon)$ is a~ball in $\mathbb{R}^D$ of radius
$\varepsilon$ centered at~$Y_i$.
\citet{smale} show that if $\lambda$ and $\varepsilon$ are chosen properly,
then $\hat M$ is homologous to $M$.
This means that $M$ and $\hat M$ share certain topological properties.
However, the result does not guarantee closeness in Hausdorff distance.
A~very relevant paper is \citet{CCDM}.
These authors consider observations generated from a~manifold and then
contaminated by additive noise as we do in Section \ref{sec::additive}.
Also, they use deconvolution methods as we do.
However, their interest is in upper bounding the Wasserstein distance
between an estimator $\hat G$ and
the distribution $G$, as a~prelude to estimating the homology of $M$.
They do not establish Hausdorff bounds.
\citet{Koltchinskii} considers estimating the number of connected components
of a~set, contaiminated by additive noise.
This corresponds to estimating the zeroth order homology.

There is a~also a~literature on estimating
principal surfaces.
A~recent paper
on this approach
with an excellent review is \citet{Principal}.
This is similar to estimating manifolds but, to the best of our knowledge,
this literature does not establish minimax bounds for estimation in
Hausdorff distance.
Finally we would like to mention the related problem
of testing for a~set of points on a~surface in a~field of uniform noise
[\citet{dots}], but, despite some similarity,
this problem is quite different.

\subsection{Notation}
We let
$B_D(x,r)$
denote a~$D$-dimensional open ball centered at $x$ with radius $r$.
If $A$ is a~set, and $x$ is a~point, then we write
$d(x,A) = \inf_{y\in A}\Vert x-y\Vert$ where
$\Vert  \cdot \Vert $ is the Euclidean norm.
Given two sets~$A$ and $B$, the \textit{Hausdorff distance} between $A$
and $B$ is
%
\begin{equation}\label{eq::hausdorff}
H(A,B) = \inf \{ \varepsilon\dvtx  A\subset B\oplus\varepsilon
   \mbox{ and }    B\subset A\oplus\varepsilon \},
\end{equation}
where
%
\begin{equation}
A\oplus\varepsilon= \bigcup_{x\in A} B_D(x,\varepsilon).
\end{equation}

The $L_1$ distance
between two distributions~$P$ and~$Q$ with densities~$p$ and~$q$
is
$\ell_1(p,q)=\int|p-q|$
and the \textit{total variation distance} between $P$ and $Q$ is
%
\begin{equation}
\mathsf{TV}(P,Q) = \sup_A~|P(A) - Q(A)|,
\end{equation}
where the supremum is over all measurable sets $A$.
Recall that
$\mathsf{TV}(P,Q) = (1/2)\ell_1(p,q)$.

Let $p(x)\wedge q(x) = \min\{p(x),q(x)\}$.
The \textit{affinity} between $P$ and $Q$ is
%
\begin{equation}
\Vert P \wedge Q\Vert  = \int p \wedge q = 1 - \frac{1}{2} \int|p-q|.
\end{equation}
Let $P^n$ denote the $n$-fold product measure
based on $n$ independent observations from $P$.
It can be shown that
%
\begin{equation}\label{eq::affinity-product}
\Vert P^n\wedge Q^n\Vert  \geq\frac{1}{8} \biggl(1 - \frac{1}{2}\int|p-q|
 \biggr)^{2n}.
\end{equation}
The convolution between two measures $P$ and $\Phi$---denoted by $P\star\Phi$---is the measure defined by
%
\begin{equation}
(P\star\Phi)(A) = \int\Phi(A-x)\, dP(x).
\end{equation}
If $\Phi$ has density $\phi$, then
$P\star\Phi$ has density
$\int\phi(y-u)\, dP(u)$.
The Fourier transform of $P$ is denoted by
%
\begin{equation}
p^* (t) = \int e^{i t^T u} \,dP(u) = \int e^{i t\cdot u} \,dP(u),
\end{equation}
where we use both $t^T u$ and $t\cdot u$ to denote the dot product.

We write $X_n = O_P(a_n)$ to mean that
for every $\varepsilon>0$, there exists $C>0$ such that
$\mathbb{P}( \Vert X_n\Vert /a_n > C) \leq\varepsilon$ for all large $n$.
Throughout, we use
symbols like
$C, C_0,C_1, c,c_0, c_1,\ldots$ to denote generic positive constants
whose value may be different in different expressions.
We write
$\mathsf{poly}(\varepsilon)$ to denote any expression of the form
$a~\varepsilon^b$ for some positive real numbers $a$ and $b$.
We write
$a_n \preceq b_n$ if
there exists $c>0$ such that
$a_n \leq c b_n$ for all large $n$.
Similarly, write
$a_n \succeq b_n$ if
$b_n \preceq a_n$.
Finally, write $a_n \asymp b_n$ if
$a_n \preceq b_n$ and
$b_n \preceq a_n$.

We will use Le Cam's lemma
to derive lower bounds,
which we now state.
This version is from
\citet{binyu}.

\begin{lemma}[(Le Cam 1973)]
\label{lemma::lecam}
Let $\cQ$
be a~set of distributions.
Let $\theta(Q)$ take values in a~metric space with metric $\rho$.
Let $Q_0,Q_1\in\cQ$ be any pair of distributions in~$\cQ$.
Let $Y_1,\ldots, Y_n$ be drawn i.i.d. from some $Q\in\cQ$
and denote the corresponding product measure by
$Q^n$.
Let $\hat\theta= \hat\theta(Y_1,\ldots, Y_n)$ be any estimator.
Then
\begin{eqnarray*}
\sup_{Q\in\cQ}
\E_{Q^n} [\rho(\hat\theta,\theta(Q)) ] & \geq&
\rho(\theta(Q_0),\theta(Q_1))   \Vert Q_0^n \wedge Q_1^n\Vert \\
& \geq&
\rho(\theta(Q_0),\theta(Q_1))   \frac{1}{8}\bigl(1- \mathsf{TV}(Q_0,Q_1)\bigr)^{2n}.
\end{eqnarray*}
\end{lemma}

\section{Assumptions}
\label{sec::conditions}

We shall be concerned with $d$-dimensional Riemannian submanifolds
of $\mathbb{R}^D$ where $d<D$.
Usually, we assume that $M$ is contained in some compact set
$\cK\subset\mathbb{R}^D$.
An exception is Section \ref{sec::additive}
where we allow noncompact manifolds.
Let $\Delta(M)$ be the largest $r$ such that
each point in $M\oplus r$ has a~unique projection onto $M$.
The quantity $\Delta(M)$ will be small
if either $M$ is not smooth or
if $M$ is close to being self-intersecting.
The quantity $\Delta(M)$
has been rediscovered many times.
It is called the \textit{condition number} in \citet{smale}
and the \textit{reach} in \citet{federer}.
Let $\cM(\kappa)$ denote all $d$-dimensional manifolds
embedded in $\mathbb{R}^D$
such that $\Delta(M) \geq\kappa$.
Throughout this paper, $\kappa$ is a~fixed positive constant.

We consider three different distributional models:

\begin{longlist}[(3)]
\item[(1)] \textit{Noiseless}.
We observe $Y_1,\ldots, Y_n \sim G$
where $G$ is supported on a~manifold $M$
where $M\in\mathcal{M} = \{ M\in\mathcal{M}(\kappa), M\subset\mathcal{K}\}$.
In this case, $Q=G$ and the observed data fall exactly on the manifold.
We assume that
$G$ has density~$g$ with respect to the uniform distribution on $M$
and that
%
\begin{equation}
0 < b(\mathcal{M}) \leq\inf_{y\in M}g(y) \leq\sup_{y\in M}g(y) \leq
B(\mathcal{M}) < \infty,
\end{equation}
where $b(\mathcal{M})$ and $B(\mathcal{M})$
are allowed to depend on the class $\mathcal{M}$, but not on the
particular manifold $M$.
Let $\cG(M)$ denote all such distributions.
In this case we define
%
\begin{equation}
\cQ= \mathcal{G} = \bigcup_{M\in\cM} \cG(M).
\end{equation}

\item[(2)] \textit{Clutter noise}.
Define $\mathcal{M}$ and $\mathcal{G}(M)$ as in the noiseless case.
We observe
%
\begin{equation}
Y_1,\ldots, Y_n \sim Q \equiv(1-\pi) U + \pi G,
\end{equation}
where
$0 < \pi\leq1$,
$U$ is uniform on the compact set
$\cK\subset\mathbb{R}^D$ and
$G\in\cG(M)$.
Define
%
\begin{equation}
\cQ=  \{ Q = (1-\pi) U + \pi G\dvtx G\in\cG(M), M \in\cM
\}.
\end{equation}

\item[(3)] \textit{Additive noise}.
In this case we allow the manifolds to be noncompact.
However, we do require that each $G$ put nontrivial probability in some
fixed compact set.
Specifically,
we again fix a~compact set $\cK$.
Let $\mathcal{M} = \mathcal{M}(\kappa)$.
Fix positive constants
$0 < b(\mathcal{M}) < B(\mathcal{M}) < \infty$.
For any $M \in\mathcal{M}$, let
$\mathcal{G}(M)$ be the set of distributions $G$ supported on $M$, such that
$G$ has density $g$ with respect to Hausdorff measure on $M$, and such that
%
\begin{equation}
0 < b(\mathcal{M}) \leq\inf_{y\in M\cap\mathcal{K}} g(y) \leq
\sup_{y\in M\cap\mathcal{K}} g(y) \leq B(\mathcal{M}) < \infty.
\end{equation}
Let $X_1, X_2, \ldots, X_n \sim G \in\cG(M)$,
and define
%
\begin{equation}\label{eq:model}
Y_i = X_i + Z_i,\qquad   i = 1, \ldots, n,
\end{equation}
where $Z_i$ are i.i.d. draws from a~distribution $\Phi$ on $\R^D$,
and where $\Phi$ is a~standard $D$-dimensional Gaussian.
Let $Q = G \star\Phi$ be the distribution of each $Y_i$
and $Q^n$ be the corresponding product measure.
Let $\cQ = \Set{G\star\Phi\dvtx  G\in\cG(M), M\in\cM}$.
\end{longlist}

These three models are an attempt
to capture the idea that we have data
falling on or near a~manifold.
These appear to be the most commonly used models.
No doubt, one could create other models as well
which is a~topic for future research.
As we mentioned earlier,
a~different noise model is considered in
\citet{smale2} and in \citet{us::2010}.
Those authors consider the case where the noise is perpendicular to the
manifold.
The former paper considers estimating the homology groups of~$M$ while
the latter paper shows that the minimax Hausdorff rate is
$n^{-{2}/{(2+d)}}$ in that case.

\section{Noiseless case}

We now derive the minimax bounds in the noiseless case.

\begin{theorem}
\label{thm::noiseless}
Under the noiseless model, we have
%
\begin{equation}
\inf_{\hat M}\sup_{Q\in\cQ}\E_{Q^n} [ H(\hat M,M)] \geq C
n^{-{2}/{d}}.
\end{equation}
\end{theorem}

\begin{pf}
Fix $\gamma>0$.
By Theorem 6 of \citet{us::2010}
there exist manifolds $M_0, M_1$
that satisfy the following conditions:
\begin{enumerate}[(3)]
\item[(1)] $M_0,M_1 \in\cM$.
\item[(2)] $H(M_0,M_1) = \gamma$.
\item[(3)] There is a~set $B\subset M_1$ such that:
\begin{enumerate}[(c)]
\item[(a)]$\inf_{y\in M_0} \norm{x-y} > \gamma/2$ for all $x\in B$.
\item[(b)]$\mu_1(B) \ge\gamma^{d/2}$ where $\mu_1$ is the uniform
measure on $M_1$.
\item[(c)] There is a~point $x\in B$ such that
$\norm{x-y} =\gamma$
where $y\in M_0$ is the closest point on $M_0$ to $x$.
Moreover, $T_x M_1$ and $T_y M_0$ are parallel
where $T_x M$ is the tangent plane to $M$ at $x$.
\end{enumerate}
\item[(4)] If $A~= \Set{y\dvtx  y\in M_1, y\notin M_0}$,
then
$\mu_1(A) \le C \gamma^{d/2}$ for some $C>0$.
\end{enumerate}
Let $Q_i = G_i$ be the uniform measure on $M_i$, for $i = 0, 1$,
and let $A$ be the set defined in the last item.
Then
$\mathsf{TV}(G_0,G_1) = G_1(A) - G_0(A) = G_1(A) \le C \gamma^{d/2}$.
From Le Cam's lemma,
%
\begin{equation}
\sup_{Q\in\cQ} \E_{Q^n} H(\hat M,M) \ge\gamma(1- \gamma^{d/2})^{2n}.
\end{equation}
Setting $\gamma= (1/n)^{2/d}$ yields the stated lower bound.
\end{pf}

\begin{figure}

\includegraphics{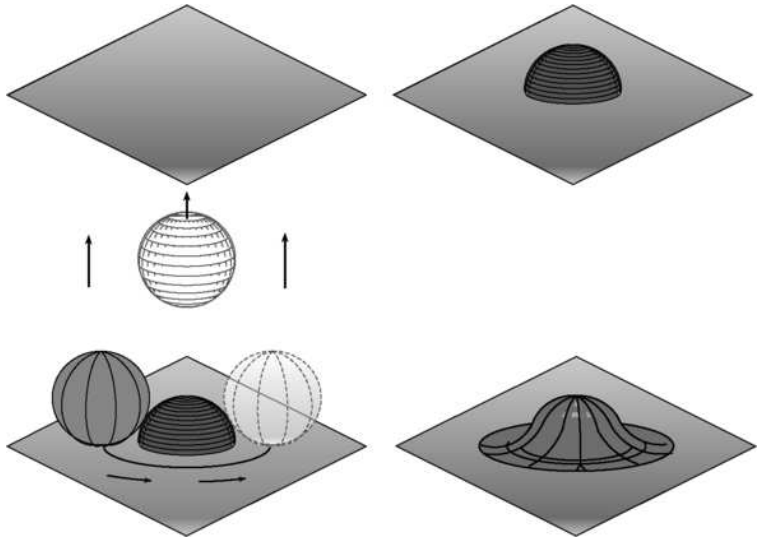}

\caption{The proof of Theorem \protect\ref{thm::noiseless}
uses two manifolds, $M_0$ and $M_1$.
A~sphere of radius $\kappa$ is pushed upward into the plane $M_0$ (top left).
The resulting manifold $M_0'$ is not smooth (top right).
A~sphere is then rolled around the manifold (bottom left) to produce a~smooth manifold $M_1$ (bottom right).
The construction is made rigorous in Theorem 6 of Genovese et~al. (\citeyear{us::2010}).}
\label{fig::flying-saucer}
\end{figure}

See Figure \ref{fig::flying-saucer}
for a~heuristic explanation of the construction of the two manifolds,
$M_0$ and $M_1$,
used in the above proof.
Now we derive an upper bound.

\begin{theorem}
\label{thm::noiseless2}
Under the noiseless model, we have
%
\begin{equation}
\inf_{\hat M}\sup_{Q\in\cQ}\E_{Q^n} [ H(\hat M,M)] \leq C
\biggl(\frac{\log n}{n} \biggr)^{{2}/{d}}.
\end{equation}
\end{theorem}

Hence, the rate is tight, up to logarithmic factors.
The proof is a~special case of the proof of the upper bound in the next
section and so is omitted.

\begin{remark*}
The Associate Editor pointed out that
the rate $(1/n)^{2/d}$ might seem counterintuitive.
For example, when $d=1$, this yields
$(1/n)^2$ which would seem to contradict the usual $1/n$ rate for
estimating the support of a~uniform distribution.
However, the slower $1/n$ rate is actually a~boundary effect
much like the boundary effects that occur
in density estimation and regression.
If we embed the uniform into $\mathbb{R}^2$
and wrap it into a~circle to eliminate the boundary,
we do indeed get a~rate of $1/n^2$.
Our assumption of smooth manifolds without boundary
removes the boundary effect.
\end{remark*}

\section{Clutter noise}

Recall that
\[
Y_1,\ldots, Y_n \sim Q=(1-\pi) U + \pi G,
\]
where
$U$ is uniform on $\mathcal{K}$,
$0 < \pi\leq1$ and $G\in\mathcal{G}$.

\begin{theorem}
\label{thm::clutter-upper}
Under the clutter model, we have
%
\begin{equation}
\inf_{\hat M}\sup_{Q\in\cQ}\E_{Q^n} [ H(\hat M,M)] \geq
C  \biggl(\frac{1}{n\pi} \biggr)^{{2}/{d}}.
\end{equation}
\end{theorem}

\begin{pf}
We define $M_0$, $M_1$ and $A$ as in the proof of Theorem
\ref{thm::noiseless}.
Let $Q_0 = (1-\pi) U + \pi G_0$ and
$Q_1 = (1-\pi) U + \pi G_1$.
Then
$\mathsf{TV}(Q_0,Q_1) = \pi\mathsf{TV}(G_0,G_1)$.
Hence
$\mathsf{TV}(Q_0,Q_1) \leq
\pi(G_1(A) - G_0(A)) = \pi G_1(A) \leq C \pi\gamma^{d/2}$.
From Le Cam's lemma,
%
\begin{equation}
\sup_{Q\in\cQ}
\E_{Q^n} [H(\hat M,M) ] \geq
\gamma(1- \pi\gamma^{d/2})^{2n}.
\end{equation}
Setting $\gamma= (1/n \pi)^{2/d}$
yields the stated lower bound.
\end{pf}

\begin{figure}

\includegraphics{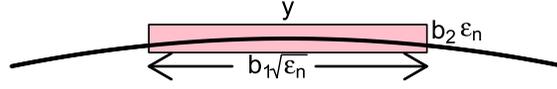}

\caption{Given a~manifold $M$ and a~point $y\in M$,
$S_M(y)$ is a~slab, centered at $y$,
with size $O(\sqrt{\varepsilon_n})$ in the $d$ directions corresponding
to the
tangent space $T_y M$ and
size $O(\varepsilon_n)$ in the $D-d$ normal directions.}
\label{fig::slab}
\end{figure}

Now we consider the upper bound.
Let $\hat Q_n$ be the empirical measure.
Let $\varepsilon_n = (K \log n/n)^{2/d}$
where $K>0$ is a~large positive constant.
Given a~manifold $M$ and a~point $y\in M$
let $S_M(y)$ denote the slab, centered at $y$,
with size $b_1 \sqrt{\varepsilon_n}$ in the $d$ directions corresponding
to the
tangent space $T_y M$ and
size $b_2 \varepsilon_n$ in the $D-d$ normal directions to the tangent space.
Here, $b_1$ and $b_2$ are small, positive constants.
See Figure \ref{fig::slab}.

Define
\[
s(M) = \inf_{y\in M} \hat Q_n[ S_M(y)] \quad  \mbox{and} \quad
\hat M_n = \mathop{\argmax}_{M} s(M).
\]
In case of ties we take any maximizer.

\begin{theorem}\label{thm::thisthis}
Let $\xi>1$ and
let $\varepsilon_n = (K \log n/n)^{2/d}$
where $K$ is a~large, positive constant.
Then
\[
\sup_{Q\in\mathcal{Q}} Q^n\bigl( H(M_0,\hat M_n) > \varepsilon_n\bigr) < n^{-\xi}
\]
and hence
\[
\sup_{Q\in\mathcal{Q}} \mathbb{E}_{Q^n}( H(M_0,\hat M_n)) \leq C
\varepsilon_n.
\]
\end{theorem}

We will use the following result, which follows
from Theorem 7 of \citet{BBL}.
This version of the result is from
\citet{SanjoyTree}.

\begin{lemma}
Let $\mathcal{A}$ be a~class of sets with VC dimension V.
Let $0 < u < 1$ and
\[
\beta_n = \sqrt{  \biggl(\frac{4}{n} \biggr)  \biggl[ V \log(2n) +
\log \biggl(\frac{8}{u} \biggr) \biggr]}.
\]
Then for all $A\in\mathcal{A}$,
\begin{eqnarray*}
&&-\min \bigl\{\beta_n \sqrt{\hat Q_n(A)}, \beta_n^2 + \beta_n
\sqrt{Q(A)} \bigr\}\\
&&\qquad  \leq
Q(A) - \hat{Q}_n(A)  \leq
\min \bigl\{ \beta_n^2 + \beta_n \sqrt{\hat Q_n(A)},  \beta_n
\sqrt{Q(A)}  \bigr\}
\end{eqnarray*}
with probability at least $1-u$.
\end{lemma}

The set of hyper-rectangles in $\mathbb{R}^D$
(which contains all the slabs)
has finite VC dimension $V$, say.
Hence, we have the following lemma obtained by setting \mbox{$u= (1/n)^\xi$}.

\begin{lemma}
Let $\mathcal{A}$ denote all hyper-rectangles in $\mathbb{R}^D$.
Let $C= 4[V + \max\{3,\xi\}]$.
Then for all $A\in\mathcal{A}$,
%
\begin{eqnarray}\label{eq::aaa}
\hat Q_n(A) & \leq& Q(A) +
\frac{C \log n}{n} + \sqrt{\frac{C \log n}{n}} \sqrt{Q(A)}
\quad\mbox{and}\\
\hat Q_n(A) & \geq& Q(A) - \sqrt{\frac{C \log n}{n}} \sqrt
{Q(A)}\label{eq::bbb}
\end{eqnarray}
with probability at least $1- (1/n)^\xi$.
\end{lemma}

Now we can prove Theorem
\ref{thm::thisthis}.

\begin{pf*}{Proof of Theorem \ref{thm::thisthis}}
Let $M_0$ denote the true manifold.
Assume that~(\ref{eq::aaa}) and (\ref{eq::bbb}) hold.
Let $y\in M_0$ and let
$A~= S_{M_0}(y)$.
Note that
$Q(A) = (1-\pi)U(A) + \pi G(A)$.
Since $y\in M_0$ and $G$ is singular, the term $U(A)$ is of lower order
and so
there exist $0 < c_1 \leq c_2 < \infty$ such that,
for all large $n$,
\[
\frac{c_1 K \log n}{n} = c_1 \varepsilon_n^{d/2} \leq
Q(A) \leq c_2 \varepsilon_n^{d/2} = \frac{c_2 K \log n}{n}.
\]
Hence
\[
\hat Q_n(A) \geq
Q(A) - \sqrt{\frac{C \log n}{n}} \sqrt{Q(A)} \geq
\frac{c_1 K \log n}{n} - \sqrt{c_2' K} \frac{\log n}{n} >
\frac{c_3 K \log n}{n}.
\]
Thus
$s(M_0) > \frac{c_3 K \log n}{n}$
with high probablity.

Now consider any $M$
for which $H(M_0,M) > \varepsilon_n$.
There exists a~point $y\in M$ such that
$d(y,M_0) > \varepsilon_n$.
It can be seen, since $M\in\mathcal{M}$, that
$S_M(y)\cap M_0 = \varnothing $.
[To see this, note that $\Delta(M) \geq\kappa>0$ implies that
the interior of any ball of radius $\kappa$ tangent to $M$ at $y$
has empty intersection with $M$ and the slab $S_M(y)$ is strictly
contained in such
a~ball for $b_1$ and $b_2$ small enough relative to $\kappa$.]
Hence
\begin{eqnarray*}
Q(S_M(y)) &=& (1-\pi)U(S_M(y)) = c_4 \varepsilon_n^{d/2} \varepsilon
_n^{D-d}\\
&=&
 \biggl(\frac{ K \log n}{n}  \biggr) c_4 \biggl( \frac{ K \log
n}{n} \biggr)^{{2(D-d)}/{d}} =
C  \biggl(\frac{\log n}{n} \biggr)^{{(2D-d)}/{d}}.
\end{eqnarray*}
So, from the previous lemma,
\begin{eqnarray*}
s(M) &=& \inf_{x\in M} \hat Q_n(S_M(x)) \leq
\hat Q_n(S_M(y))\\
& \leq &
Q(S_M(y)) + \frac{C \log n}{n} + \sqrt{\frac{C \log n}{n}} \sqrt
{Q(S_M(y))}\\
&=&
 \biggl(\frac{K \log n}{n} \biggr)^{{(2D-d)}/{d}} +
\frac{C \log n}{n} +  \biggl(\frac{K \log n}{n} \biggr)^{{D}/{d}}
< \frac{C_3 K \log n}{n} = s(M_0)
\end{eqnarray*}
since $D>d$ and $K$ is large.
Let $\mathcal{M}_n = \{M\in\mathcal{M}\dvtx H(M_0,M) > \varepsilon_n\}$.
We conclude that
\[
Q^n \bigl( s(M) > s(M_0) \mbox{ for some }M\in\mathcal{M}_n \bigr) <
 \biggl(\frac{1}{n} \biggr)^{\xi}.
\]
\upqed\end{pf*}

\section{Additive noise}
\label{sec::additive}

Let us recall the model.
Let $\mathcal{M} = \mathcal{M}(\kappa)$.
We allow the manifolds to be noncompact.
Fix positive constants
$0 < b(\mathcal{M}) < B(\mathcal{M}) < \infty$.
For any $M \in\mathcal{M}$ let
$\mathcal{G}(M)$ be the set of distributions $G$ supported on $M$ such that
$G$ has density $g$ with respect to Hausdorff measure on $M$ and such that
%
\begin{equation}
0 < b(\mathcal{M}) \leq\inf_{y\in M\cap\mathcal{K}} g(y) \leq
\sup_{y\in M\cap\mathcal{K}} g(y) \leq B(\mathcal{M}) < \infty,
\end{equation}
where $\cK$
is a~compact set.
Let $X_1, X_2, \ldots, X_n \sim G \in\cG(M)$,
and define
%
\begin{equation}\label{eq:model2}
Y_i = X_i + Z_i, \qquad  i = 1, \ldots, n,
\end{equation}
where $Z_i$ are i.i.d. draws from a~distribution $\Phi$ on $\R^D$,
and where $\Phi$ is a~standard $D$-dimensional Gaussian.
Let $Q = G \star\Phi$ be the distribution of each $Y_i$
and $Q^n$ be the corresponding product measure.
Let $\cQ = \Set{G\star\Phi\dvtx  G\in\cG(M), M\in\cM}$.

Since we allow the manifolds to be noncompact,
the Hausdorff distance could be unbounded.
Hence we define a~truncated loss function,
%
\begin{equation}
\label{eq:loss}
L(M, \hat M) = H(M \intersect\cK, \hat M \intersect\cK).
\end{equation}

\begin{theorem}
\label{thm::additive-minimax}
For all large enough $n$,
%
\begin{equation}
\inf_{\hat M}\sup_{Q\in\cQ}\E_Q [L(M,\hat M)] \ge\frac{C}{\log n}.\vadjust{\goodbreak}
\end{equation}
\end{theorem}

\begin{pf}
Define $\tilde c\dvtx \R\to\R$ and $c\dvtx \R^d \to\R^{D-d}$ as follows:
$\tilde c(x) =\break  \cos(x/(a\sqrt{\gamma}))$ and
$c(u) =  (\prod_{\ell=1}^d \tilde c(u_\ell), 0, \ldots,
0 )^T$.
Let
$M_0 = \Set{ (u,\gamma  c(u))\dvtx  u\in\R^d}$
and
$M_1 = \Set{ (u,-\gamma  c(u))\dvtx u\in\R^d}.$
See Figure \ref{fig::two-cosines} for a~picture of~$M_0$ and~$M_1$ when $D = 2, d = 1$.
Later, we will show that $M_0,M_1 \in\cM$.

\begin{figure}

\includegraphics{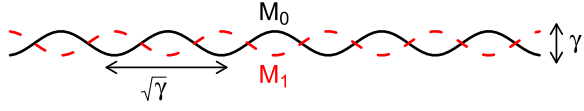}

\caption{The two least favorable manifolds $M_0$ and $M_1$
in the proof of Theorem \protect\ref{thm::additive-minimax}
in the special case where $D=2$ and $d=1$.}
\label{fig::two-cosines}
\end{figure}

Let $U$ be a~$d$-dimensional random variable with density $\zeta$
where $\zeta$ is $d$-dimensional standard Gaussian density.
Let $\tilde\zeta$ be a~one-dimensional $N(0,1)$ density.
And define $G_0$ and $G_1$
by $G_0(A) = \P((U,\gamma c(U))\in A)$ and
$G_1(A) = \P((U,-\gamma c(U))\in A)$.

We begin by bounding
$\int|q_1 - q_0|^2$.
Define the $D$-cube $\mathcal{Z} = [-1/(2a\sqrt{\gamma}),\break 1/(2a\sqrt
{\gamma})]^D$.
Then, by Parseval's identity,
and that fact that $q_j^* = \phi^* g_j^*$,
\begin{eqnarray*}
(2\pi)^D \int|q_1 - q_0|^2 &=& \int|q_1^* - q_0^*|^2 = \int|\phi
^*|^2  |g_1^* - g_0^*|^2 \\
&=& \int_{\mathcal{Z}} |\phi^*|^2  |g_1^* - g_0^*|^2 + \int_{\mathcal{Z}^c} |\phi^*|^2  |g_1^* - g_0^*|^2 \\
&\equiv&        I               +             \mathit{II}.
\end{eqnarray*}
Then
\begin{eqnarray*}
\mathit{II}
&=& \int_{\mathcal{Z}^c} |g_1^*(t) - g_0^*(t)|^2   |\phi^*(t)|^2\\
&\le&\int_{\mathcal{Z}^c} |\phi^*(t)|^2
\le C  \biggl(\int_{1/(2a\sqrt{\gamma})}^\infty e^{-t^2} \,d t
\biggr)^D \\
&\le&\mathsf{poly}(\gamma)  e^{-D/4 a^2\gamma}.
\end{eqnarray*}

Now we bound $I$.
Write $t\in\R^D$ as $(t_1,t_2)$
where
$t_1 = (t_{11}, \ldots, t_{1d}) \in\mathbb{R}^d$ and
$t_2 = (t_{21}, \ldots, t_{2(D-d)}) \in\mathbb{R}^{D-d}$.
Let
$c_1(u) = \prod_{\ell=1}^d \tilde{c}(u_\ell)$
denote the first component of the vector-valued function $c$.
We have
\begin{eqnarray*}
g_1^*(t) - g_0^*(t)
&=& \int_{\mathbb{R}^d}  \bigl( e^{i t_1\cdot u + i t_{21} \gamma
c_1(u)}-e^{i t_1\cdot u - i t_{21} \gamma c_1(u)}  \bigr) \zeta(u)\,du\\
&=& 2i \int e^{i t_1\cdot u} \sin( t_{21} \gamma c_1(u)) \zeta(u) \,du\\
&= &2i \int e^{i t_1 \cdot u}
\sum_{k=0}^\infty\frac{ (-1)^k t_{21}^{2k+1} \gamma^{2k+1}}{
(2k+1)!} c_1^{2k+1}(u) \zeta(u) \,du\\
&=& 2i \sum_{k=0}^\infty\frac{ (-1)^k t_{21}^{2k+1} \gamma^{2k+1}}{
(2k+1)!} \int e^{i t_1\cdot u} c_1^{2k+1}(u) \zeta(u) \,du\\
&=& 2i \sum_{k=0}^\infty\frac{ (-1)^k t_{21}^{2k+1} \gamma^{2k+1}}{ (2k+1)!}
\prod_{\ell=1}^d \int e^{i t_{1\ell} u_\ell} {\tilde
c}^{2k+1}(u_\ell) \tilde\zeta(u_\ell) \,du_\ell\\
&=& 2i \sum_{k=0}^\infty\frac{ (-1)^k t_{21}^{2k+1} \gamma^{2k+1}}{
(2k+1)!} \prod_{\ell=1}^d ({\tilde c}^{2k+1}\tilde\zeta)^*
(t_{1\ell})\\
&=& 2i \sum_{k=0}^\infty\frac{ (-1)^k t_{21}^{2k+1} \gamma^{2k+1}}{
(2k+1)!} \prod_{\ell=1}^d m_k(t_{1\ell}),
\end{eqnarray*}
where
%
\begin{equation}
m_k(t_{1\ell})= ({\tilde c}^{2k+1} \tilde\zeta)^*(t_{1\ell}) =
 (\underbrace{{\tilde c}^*\star{\tilde c}^* \star
\cdots\star{\tilde c}^*}_{2k+1 \ \mathrm{times}}\star{\tilde\zeta
}^* )(t_{1\ell}).
\end{equation}
Note that
\[
{\tilde c}^* = \tfrac{1}{2}\delta_{{-1}/{(a\sqrt{\gamma})}} +
\tfrac{1}{2}\delta_{1/{(a\sqrt{\gamma})}},
\]
where $\delta_y$ a~Dirac delta function at $y$, that is, a~generalized
function corresponding to point evaluation at $y$.
For any integer $r$, if we
convolve ${\tilde c}^*$ with itself $r$ times, we have that
%
\begin{equation}
\underbrace{{\tilde c}^*\star{\tilde c}^* \star\cdots\star{\tilde
c}^*}_{r\ \mathrm{times}} =
 \biggl(\frac{1}{2} \biggr)^r \sum_{j=0}^r \pmatrix{r\cr j} \delta_{a_j},
\end{equation}
where $a_j = (2j-r)/(a\sqrt{\gamma})$.
Thus
%
\begin{equation}
m_k(t_{1\ell}) =  \biggl(\frac{1}{2} \biggr)^{2k+1} \sum_{j=0}^{2k+1}
\pmatrix{2k+1\cr j} {\tilde\zeta}^*(t_{1\ell} - a_j).
\end{equation}
Now
${\tilde\zeta}^* (t_{1\ell}) = \exp ( - \frac{t_{1\ell
}^2}{2}  )$
and ${\tilde\zeta}^*(s) \le1$ for all $s\in\R$.
For $t\in\mathcal{Z}$,
${\tilde\zeta}^*(t_{1\ell} - a_j) \le e^{-{1}/{(2a^2\gamma)}}$,
and thus
$|m_k(t_{1\ell})| \le e^{-{1}/{(2a^2\gamma)}}$.
Hence,
$\prod_{\ell=1}^d |m_k(t_{1\ell})| \le e^{-{d}/{(2a^2\gamma)}}$.
It follows that for $t\in\mathcal{Z}$,
\begin{eqnarray*}
|g_1^*(t) - g_0^*(t)|
&\le& 2 \sum_{k=0}^\infty\frac{ |t_{21}|^{2k+1} \gamma^{2k+1}}{
(2k+1)!} \prod_{\ell=1}^d |m_k(t_{1\ell})|\\
&\le& e^{-{d}/{(2a^2\gamma)}} \sum_{k=0}^\infty\frac{
|t_{21}|^{2k+1} \gamma^{2k+1}}{ (2k+1)!}\\
&\le& e^{-{d}/{(2a^2\gamma)}}  \operatorname{sinh}(|t_{21}|\gamma) \le
e^{-{d}/{(2a^2\gamma)}}.
\end{eqnarray*}
So,
\begin{eqnarray*}
I &=& \int_{\mathcal{Z}} |g_1^*(t) - g_0^*(t)|^2 |\phi^*(t)|^2 \, d
t\\
&\leq& \int_{\mathcal{Z}} |g_1^*(t) - g_0^*(t)|^2 \, d  t\\
&\le&\mathsf{Volume}(\mathcal{Z}) e^{-{d}/{(a^2\gamma)}} = \mathsf{
poly}(\gamma)e^{-{d}/{(a^2\gamma)}}.
\end{eqnarray*}
Hence,
%
\begin{eqnarray}\nonumber
\int|q_1 - q_0|^2 & \le& I + \mathit{II} \leq
\mathsf{poly}(\gamma) e^{-d/a^2\gamma} +
\mathsf{poly}(\gamma) e^{-D/4a^2\gamma}\\
& = & \mathsf{poly}(\gamma) e^{- 2w/\gamma},
\end{eqnarray}
where
$2w = \min\{ d/a^2, D/(4 a^2)\}$.

Next we bound
$\int|q_1 - q_0|$ so that we can apply Le Cam's lemma.
Let $T_\gamma$ be a~ball centered at the origin with radius $1/\gamma$.
Then, by Cauchy--Schwarz,
\begin{eqnarray*}
\int|q_1 - q_0| &=&
\int_{T_\gamma} |q_1 - q_0| + \int_{T_\gamma^c} |q_1 - q_0|\\
& \leq& \sqrt{\mathsf{Volume}(T_\gamma)} \sqrt{ \int|q_1-q_0|^2} +
\int_{T_\gamma^c} |q_1 - q_0|\\
& \leq& \mathsf{poly}(\gamma) e^{-{w}/{\gamma}} + \int
_{T^c_\gamma} |q_1 - q_0|.
\end{eqnarray*}
For all small $\gamma$ we have that
$\mathcal{K} \subset T_\gamma$.
Hence,
\begin{eqnarray*}
\int_{T_\gamma^c} |q_1 - q_0| &\leq&
\int_{M_1}\int_{T_\gamma^c} \phi(\Vert y-u\Vert ) + \int_{M_0}\int
_{T_\gamma^c} \phi(\Vert y-u\Vert ) \leq
\mathsf{poly}(\gamma)e^{-{D}/{\gamma^2}}\\
& \leq& \mathsf{poly}(\gamma)e^{-{w}/{\gamma}}.
\end{eqnarray*}
Putting this all together, we have that
$\int|q_1 - q_0| \leq
\mathsf{poly}(\gamma) e^{-{w}/{\gamma}}$.

Now we apply Lemma \ref{lemma::lecam}
and conclude that, for every $\gamma>0$,
\[
\sup_Q \E(L(M,\hat M)) \ge\frac{\gamma}{8} \bigl(1- \mathsf{poly}(\gamma)
e^{-w/\gamma}\bigr)^{2n}.
\]
Set $\gamma\asymp w/\log n$ and conclude that, for all large $n$,
\[
\sup_Q \E(L(M,\hat M)) \ge\frac{w}{8 e^2} \frac{1}{\log n}.
\]
This concludes the
proof of the lower bound
except that
it remains to show
that $M_0, M_1 \in\cM(\kappa)$.
Note that
$|{\tilde c}''(u)| = a^{-2}|\cos(u/(a\sqrt{\gamma})|$.
Hence, as long as \mbox{$a~> \sqrt{\kappa}$},
$\sup_u |{\tilde c}''(u)| < 1/\kappa$.
It now follows that $M_0,M_1\in\cM(\kappa)$. 
This completes the proof.
\end{pf}

\begin{remark*}
Consider the special case where
$D=2$, $d=1$ and the manifold
has the special form
$\{(u,m(u))\dvtx u\in\mathbb{R}\}$ for some function
$m\dvtx \mathbb{R}\to\mathbb{R}$.
In this case, estimating the manifold
is like estimating a~regression function with errors in variables.
(More on this in Section \ref{section::deconv}.)
The rate obtained for
estimating a~regression function with errors in variables
under these conditions
[\citet{FanTruong}] is $1/\log n$
in agreement with our rate.
However, the proof technique is not quite the same
as we explain in Section~\ref{section::deconv}.
\end{remark*}

\begin{remark*}
The proof of the lower bound is similar to
other lower bounds in deconvolution problems.
There is an interesting technical difference, however.
In standard deconvolution, we can choose $G_0$ and $G_1$
so that $g_1^*(t)-g_0^*(t)$ is zero
in a~large neighborhood around the origin.
This simplifies the proof considerably.
It appears we cannot do this in the manifold case since~$G_0$ and~$G_1$ have different supports.
\end{remark*}

Next we construct an upper bound.
We use a~standard deconvolution density estimator $\hat g$
(even thought $G$ has no density),
and then we threshold this estimator.

\begin{theorem} \label{theorem::upper-bound}
Fix any $0 < \delta< 1/2$.
Let $h = 1/\sqrt{\log n}$.
Let $\lambda_n$ be such that
\[
C'  \biggl(\frac{1}{h} \biggr)^{D-d} < \lambda_n < C''  \biggl(\frac
{1}{L} \biggr)^{2k}
 \biggl(\frac{1}{h} \biggr)^{D-d},
\]
where
$k \geq d/(2\delta)$,
$C'$ is defined in
Lemma \ref{lemma::inf-over-M}
and $C''$ and $L$ are defined in Lemma~\ref{lemma::sup-off-M}.
Define
$\hat M = \{y\dvtx \hat g(y) > \lambda_n\}$
where $\hat g$ is defined in (\ref{eq::hat-g}).
Then for all large $n$,
%
\begin{equation}
\inf_{\hat M}\sup_{Q\in\cQ}\E_Q [ L(M,\hat M)] \le
C  \biggl(\frac{1}{\log n} \biggr)^{{(1-\delta)}/{2}}.
\end{equation}
\end{theorem}

Let us now define the estimator in more detail.
Define
$\psi_k(y)=\mathrm{sinc}^{2k}(y/(2k))$.
By elementary calculations, it follows that
\[
\psi_k^*(t) = 2k B_{2k} \biggl(\frac{t}{2k} \biggr),
\]
where
$B_r = \underbrace{J\star\cdots\star J}_{r\  \mathrm{times}}$
where $J=\frac{1}{2}I_{[-1,1]}$.
The following properties of~$\psi_k$ and~$\psi_k^*$
follow easily:
\begin{enumerate}[(3)]
\item[(1)] The support of $\psi_k^*$ is $[-1,1]$.
\item[(2)] $\psi_k \geq0$ and $\psi_k^* \geq0$.
\item[(3)] $\int\psi_k^*(t)\, dt = \psi_k(0) = 1$.
\item[(4)] $\psi_k^*$ and $\psi_k$ are spherically symmetric.
\item[(5)] $|\psi_k(y)| \leq 1/( (2k)^{2k}|y|^{2k})$ for all $|y| >
\pi/(2k)$.
\end{enumerate}
Abusing notation somewhat, when $u$ is a~vector, we take $\psi_k(u)
\equiv\psi_k(\norm{u})$.\vadjust{\goodbreak}

Define
%
\begin{equation}
\hat g^*(t) = \frac{\hat q^*(t)}{\phi^*(t)}\psi^*_k(h t),
\end{equation}
where
$\hat q^*(t) = \frac{1}{n}\sum_{i=1}^n e^{-i t^T Y_i}$
is the empirical characteristic function.
Now define
%
\begin{equation}\label{eq::hat-g}
\hat g(y) =  \biggl(\frac{1}{2 \pi} \biggr)^D
\int e^{-i t^T y} \frac{\psi^*_k(h t)\hat q^*(t)}{\phi^*(t)}\, dt.
\end{equation}
Let
$\overline{g}(y) = \E(\hat g(y))$.

\begin{lemma}
For all $y\in\mathbb{R}^D$,
\[
\overline{g}(y) =  \biggl(\frac{1}{2\pi h} \biggr)^D \int\psi_k
\biggl( \frac{\Vert y-u\Vert }{h} \biggr)\, dG(u).
\]
\end{lemma}

\begin{pf}
Let $\psi_{k,h}(x) = h^{-D}\psi_k(x/h)$.
Hence,
$\psi^*_{k,h}(t) = \psi_{k}^*(th)$.
Now,
\begin{eqnarray*}
\overline{g}(y) &=&
 \biggl(\frac{1}{2\pi} \biggr)^D\int e^{-i t^T y} \frac{\psi^*_k(th)
q^*(t)}{\phi^*(t)}\, dt\\[-2pt]
&=&
 \biggl(\frac{1}{2\pi} \biggr)^D\int e^{-i t^T y} \frac{\psi
^*_k(th)g^*(t)\phi^*(t) }{\phi^*(t)} \,dt\\[-2pt]
&=&
 \biggl(\frac{1}{2\pi} \biggr)^D\int e^{-i t^T y} \psi^*_k(th)g^*(t)
\,dt\\[-2pt]
&=&  \biggl(\frac{1}{2\pi} \biggr)^D\int e^{-i t^T y} \psi
^*_{k,h}(t)g^*(t)\, dt=
 \biggl(\frac{1}{2\pi} \biggr)^D\int e^{-i t^T y} (g\star\psi
_{k,h})^*(t) \,dt\\[-2pt]
&=&  \biggl(\frac{1}{2\pi} \biggr)^D(g\star\psi_{k,h})(y)=
 \biggl(\frac{1}{2\pi} \biggr)^D\int\psi_{k,h}(y-u) \,dG(u)\\[-2pt]
&=& \frac{1}{h^D} \biggl(\frac{1}{2\pi} \biggr)^D \int\psi_k \biggl(
\frac{y-u}{h} \biggr)\, dG(u).
\end{eqnarray*}
\upqed\end{pf}

\begin{lemma}
\label{lemma::inf-over-M}
We have that
$\inf_{y\in M\intersect\cK}\overline{g}(y) \ge C' h^{d - D}$.
\end{lemma}

\begin{pf}
Choose any $x\in M\intersect\cK$ and let $B = B(x, C h)$.
Note that $G(B) \geq b(\mathcal{M}) c h^d$.
Hence,
\begin{eqnarray*}
{\overline g}(x) & = &
(2\pi)^{-D}
h^{-D} \int\psi_k \biggl( \frac{x - u}{h}  \biggr) \, d G(u)\\[-2pt]
& \geq&
(2\pi)^{-D}h^{-D} \int_B \psi_k \biggl( \frac{x - u}{h}  \biggr) \,d
G(u)\\[-2pt]
& \geq&
(2\pi)^{-D} h^{-D} G(B) = C' h^{d - D}.
\end{eqnarray*}
\upqed\end{pf}\eject

\begin{lemma}
\label{lemma::sup-off-M}
Fix $0 < \delta< 1/2$.
Suppose that $k \ge d/(2\delta)$.
Then,
%
\begin{equation}
\sup\Set{\overline{g}(y)\dvtx  y\in\cK,  d(y,M) > L h^{1-\delta}}
\le
C'' L^{-2k} \biggl(\frac{1}{h} \biggr)^{D-d}.
\end{equation}
\end{lemma}

\begin{pf}
Let $y$ be such that
$d(y,M) > L h^{1-\delta}$.
For integer $j \ge1$, define
\[
A_j = \bigl[B\bigl(y,(j+1)L h^{1-\delta}\bigr)- B(y,j L h^{1-\delta})\bigr]\intersect M
\intersect\cK.
\]
Then
\begin{eqnarray*}
\overline{g}(y)
&=&  \biggl(\frac{1}{2\pi h} \biggr)^D \int\psi_k \biggl( \frac
{\Vert u-y\Vert }{h} \biggr) \,dG(u)\nonumber\\
&\le&  \biggl(\frac{1}{2\pi h} \biggr)^D \sum_{j=1}^\infty\int_{A_j}
\psi_k \biggl( \frac{\Vert u-y\Vert }{h} \biggr) \,dG(u)\nonumber\\
&\le& \biggl(\frac{1}{2\pi h} \biggr)^D \sum_j \int_{A_j}  \biggl(
\frac{2k h}{\Vert u-y\Vert } \biggr)^{2k} \,dG(u)\nonumber\\
&\le& C \biggl(\frac{1}{h} \biggr)^D \sum_j \int_{A_j}  \biggl( \frac
{h}{j L h^{1-\delta}} \biggr)^{2k} \,dG(u)\nonumber\\
&\le& C \biggl(\frac{1}{h} \biggr)^D L^{-2k}   h^{2k\delta}  \sum_j
 \biggl( \frac{1}{j} \biggr)^{2k} G(A_j) \nonumber\\\label{*}
\hspace*{-84pt}\mbox{(*)}\hspace*{84pt}&\le& C \biggl(\frac{1}{h} \biggr)^D L^{-2k}   h^{2k\delta} \\
\label{**}
\hspace*{-84pt}\mbox{(**)}\hspace*{78.5pt}&\le& C \biggl(\frac{1}{h} \biggr)^D L^{-2k}   h^{d} \\
&\le& C'' L^{-2k}  \biggl(\frac{1}{h} \biggr)^{D-d},\nonumber
\end{eqnarray*}
where equation (\hyperref[*]{*}) follows because $G$ is a~probability measure
and $\sum_j j^{-2k}< \infty$,
and equation (\hyperref[**]{**}) follows because $2k \delta\ge d$.
\end{pf}

Now define
$\Gamma_n = \sup_y |\hat{g}(y) - \overline{g}(y)|$.

\begin{lemma}
Let $h = 1/\sqrt{\log n}$, and let $\xi> 1$.
Then, for large $n$,
%
\begin{equation}
\Gamma_n =  \biggl(\frac{1}{\sqrt{\log n}} \biggr)^{4k + 4-D}
\end{equation}
on an event $\cA_{n}$ of probability at least $1-n^{-\xi}$.
\end{lemma}

\begin{pf}
We proceed as in Theorem 2.3 of
\citet{Stefanski1990229}.
Note that
%
\begin{equation}
\hat g(y) - \overline{g}(y) =
 \biggl(\frac{1}{2\pi} \biggr)^D\int e^{-i t^T y} \frac{\psi
_k^*(th)}{\phi^*(t)} \bigl(\hat q^*(t) - q^*(t)\bigr) \,dt,
\end{equation}
and also note that the integrand is 0 for $\Vert t\Vert  > 1/h$.
So
%
\begin{equation}\label{eq::sup-bound}
\sup_y|\hat g(y) - \overline{g}(y)| \le
\frac{\Delta_n}{(2\pi)^D}  \biggl|\int_{\Vert t\Vert \leq1/h} \frac{\psi
_k^*(th)}{\phi^*(t)}\,dt  \biggr|,
\end{equation}
where
$\Delta_n =\sup_{\Vert t\Vert  < 1/h} |\hat q^*(t) - q^*(t)|$.

For $D = 1$, it follows from Theorem 4.3 of
\citet{Yukich1985245} that
%
\begin{equation}
Q^n(\Delta_n > 4\varepsilon) \le
4 N(\varepsilon) \exp \biggl( - \frac{n \varepsilon^2}{8 + 4\varepsilon
/3} \biggr) +
8 N(\varepsilon) \exp \biggl( - \frac{n \varepsilon}{96} \biggr),
\end{equation}
where $N(\varepsilon)$ is the bracketing number of the set of complex
exponentials, which\vspace*{1pt} is given by
$N(\varepsilon) = 1 + \frac{24 M_\varepsilon T_n}{\varepsilon}$,
and $M_\varepsilon$ is defined by $Q(\norm{Y} > M_\varepsilon) \le
\varepsilon/4$.
By a~similar argument, we have that in $D$ dimensions,
%
\begin{equation}
\sup_{Q\in\cQ_n} Q^n(\Delta_n > 4\varepsilon) \le
4 N(\varepsilon) \exp \biggl( - \frac{n \varepsilon^2}{8 + 4\varepsilon
/3} \biggr) +
8 N(\varepsilon) \exp \biggl( - \frac{n \varepsilon}{96} \biggr),
\end{equation}
where now
%
\begin{equation}
N(\varepsilon) = C  \biggl[1 + \frac{24 M_\varepsilon T_n}{\varepsilon}
\biggr] \varepsilon^{-(D-1)},
\end{equation}
and $M_\varepsilon$ is defined by $\sup_{Q\in\mathcal{Q}_n} Q(\norm{Y} >
M_\varepsilon) \le\varepsilon/4$.
Note that
$M_\varepsilon=O(1)$.
It follows that $\Delta_n \le\sqrt{\frac{C\log n}{n}}$
except on a~set of probability $n^{-\xi}$ where $\xi$ can
be made arbitrarily large by taking $C$ large.

Now, note that ${\psi_k^*(h t)}/{\phi^*(t)}$ is a~spherically symmetric
function $R(\norm{t})$.
Hence,
\[
\int_{\norm{t}\le1/h}  \frac{\psi_k^*(h t)}{\phi^*(t)} \,d t =
C \int_{s=0}^{1/h} R(s) s^{D - 1}\, d s \leq
C h^{4k + 4 - D} e^{{1}/{(2 h^2)}},
\]
where the last result follows from Lemma 3.1 in Stefanski (\citeyear{Stefanski1990229})
using parameters $\delta=2$, $\gamma=1/2$, $r = 2k + 2$,
$\beta= D - 1$, with $\lambda=h$.
The value of $r$ follows from the definition of $\psi_k^*$.
The result now follows by
combining this bound with (\ref{eq::sup-bound}).
\end{pf}

Now we can complete the proof of the upper bound.

\begin{pf*}{Proof of Theorem \protect\ref{theorem::upper-bound}}
On the event $\cA_n$
where $\Gamma_n \leq(1/\sqrt{\log n})^{4k+4-D}$
(defined in the previous lemma),
we have
\begin{eqnarray*}
\inf_{y\in M\intersect\cK} \hat g(y) &\geq&
\inf_{y\in M\intersect\cK} \overline{g}(y) - \Gamma_n \geq
C  \biggl(\frac1h \biggr)^{D-d} -  \biggl(\frac1{\sqrt{\log n}}
\biggr)^{4k + 4 - D} \\
& \geq& (C/2)  \biggl(\frac1h \biggr)^{D-d} > \lambda_n.
\end{eqnarray*}
This implies that $M\intersect\cK\subset\hat M \intersect\cK$

Next, we have
\begin{eqnarray*}
\mathop{\sup_{y\in\cK}}_{d(y,M)\ge L h^{1-\delta}}
 \hat g(y)
&\le&
\mathop{\sup_{y\in\cK}}_{d(y,M)\ge L h^{1-\delta}} \overline
{g}(y) + \Gamma_n \\
&\leq&
C L^{-2k}  \biggl(\frac1h \biggr)^{D-d} +  \biggl(\frac1{\sqrt{\log
n}} \biggr)^{4k + 4 - D} \\
& \leq& 2 C L^{-2k}  \biggl(\frac1h \biggr)^{D-d} < \lambda_n
\end{eqnarray*}
for large enough $L$.
This implies that
\[
\Set{y\dvtx  y\in\cK \mbox{ and } d(y,M) \ge L h^{1-\delta}} \cap
\hat M = \varnothing .
\]

Therefore, on $\cA_n$,
$L(M, \hat M) \le C  (\frac{1}{\log n} )^{{(1-\delta)}/{2}}$
and hence,
\begin{eqnarray*}
\E(L(M,\hat M)) &=& \E(L(M,\hat M) 1_{\cA_n}) + \E(L(M,\hat M)
1_{\cA_n^c}) \\
& \leq&
C  \biggl(\frac{1}{\log n} \biggr)^{{(1-\delta)}/{2}} + Q^n(\cA
_n^c) \\
& \leq&
C  \biggl(\frac{1}{\log n} \biggr)^{{(1-\delta)}/{2}} + n^{-\xi}
\leq
C  \biggl(\frac{1}{\log n} \biggr)^{{(1-\delta)}/{2}},
\end{eqnarray*}
and the theorem is proved.
\end{pf*}

\begin{remark*}
Again, the proof of the upper bound is similar to
proofs used in other deconvolution problems.
But once more, there are interesting differences.
In particular, the density estimator $\hat g$ is not estimating
any underlying density since the measure $G$ is singular and thus
does not have a~density.
Hence, the usual bias calculation is meaningless.
\end{remark*}

\begin{remark*}
Note that
$\hat M$ is a~set not a~manifold;
if desired, we can replace $\hat M$ with any
manifold
in $\{M\in\mathcal{M}\dvtx M \subset\hat M\}$,
and then the estimator is a~manifold and the rate is the same.
\end{remark*}

\begin{remark*}
The upper bound is slightly slower than the lower bound.
The rate is consistent with the results
in \citet{CCDM} who show that
$\mathbb{E}(W_2(\hat g,G)) \leq C /\sqrt{\log n}$
where $W_2$ is the Wasserstein distance.
In the special case where
the manifold
has the form
$\{(u,m(u))\dvtx u\in\mathbb{R}\}$ for some function~$m$,
the problem can be viewed as nonparametric regression with measurement
error; see Section \ref{section::deconv}.\vadjust{\goodbreak}
In this special case,
we can use the deconvolution
kernel regression estimator in \citet{FanTruong}
which achieves the rate $1/\log n$.
We do not know of any estimator in the general case that achieves the rate
$1/\log n$, although we conjecture that the following estimator
might have a~better rate: let $(\hat M, \hat G)$ minimize
$\sup_{\Vert t\Vert \leq T_n} |\hat q^*(y) - q^*_{M,G}(t)|$ where $T_n =
O(\sqrt{\log n})$.
In any case, as with all Gaussian deconvolution problems, the rate is
very slow,
and the difference between $1/\log n$ and $1/\sqrt{\log n}$ is not of
practical consequence.
\end{remark*}

\section{Singular deconvolution}
\label{section::deconv}

Estimating a~manifold under additive noise
is related to deconvolution.
It is also related to
regression with errors in variables.
The purpose of this section is to explain the connections
between the problems.

\subsection{Relationship to density deconvolution}
Recall that the model is
$Y= X + Z$
where $X\sim G$,
$G$ is supported on a~manifold $M$ and $Z \sim\Phi$.
$G$ is a~singular measure supported on the $d$-dimensional manifold $M$.

Now consider a~somewhat simpler model:
suppose again that $Y_i = X_i + Z_i$,
but suppose that $X$ has a~density $g$ on $\mathbb{R}^D$
(instead of being supported on a~manifold).
All three distributions
$Q$, $G$ and $\Phi$ have $D$-dimensional support
and $Q = G\star\Phi$.
The problem of recovering the density $g$ of $X$
from $Y_1,\ldots, Y_n$ is the usual
density deconvolution problem.
A~key reference is \citet{Fan}.

Most of the existing literature
on deconvolution assumes that $X$ and $Y$ have the same support,
or at least that the supports have the same dimension; an exception is
\citet{Koltchinskii}.
Manifold learning may be regarded as the problem of
deconvolution for singular measures.

\begin{figure}

\includegraphics{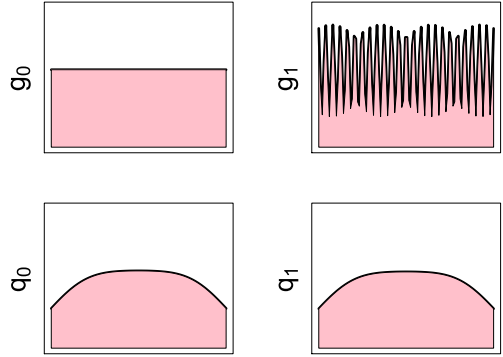}%
\vspace*{-3pt}
\caption{A~typical least favorable pair
for proving a~lower bounds in ordinary deconvolution.
The top left plot is a~density $g_0$ and
the top right plot is a~density $g_1$
which is a~perturbed version of $g_0$.
The $L_1$ distance between the densities is $\varepsilon$.
The bottom plots are
$q_0 = \int\phi(y-x) g_0(x) \,dx$ and
$q_1 = \int\phi(y-x) g_1(x)\, dx$.
These densities are nearly indistinguishable and, in fact,
their total variation distance is $e^{-1/\varepsilon}$.}
\label{fig::picture2}\vspace*{-3pt}
\end{figure}

It is instructive to compare the least favorable pair
used for proving the lower bounds
in the ordinary case versus the singular case.
Figure \ref{fig::picture2}
shows a~typical least favorable pair
for proving a~lower bound in ordinary deconvolution.
The top left plot is a~density $g_0$, and
the top right plot is a~density $g_1$
which is a~perturbed version of $g_0$.
The $L_1$ distance between the densities is~$\varepsilon$.
The bottom plots are
$q_0 = \int\phi(y-x) g_0(x) \,dx$ and
$q_1 = \int\phi(y-x) g_1(x) \,dx$.
These densities are nearly indistinguishable, and, in fact,
their total variation distance is of order $e^{-1/\varepsilon}$.
Of course, these distributions have the same support
and hence such a~least favorable pair will not suffice
for proving lower bounds in the manifold case
where we will need two densities with different support.

\begin{figure}[b]
\vspace*{-3pt}
\includegraphics{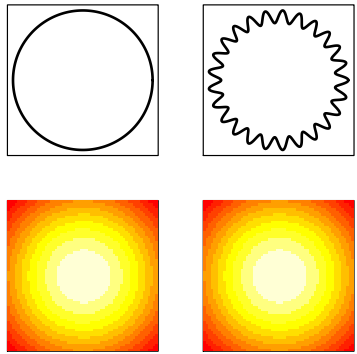}%
\vspace*{-3pt}
\caption{
The type of least favorable pair
needed for proving lower bounds in manifold learning.
The distribution $g_0$ is uniform on the circle in the top left plot.
The distribution~$g_1$ is uniform on the perturbed circle in the top
right plot.
The Hausdorff distance between the supports of the densities is~$\varepsilon$.
The bottom plots are heat maps of
$q_0 = \int\phi(y-x) g_0(x) \,dx$ and
$q_1 = \int\phi(y-x) g_1(x) \,dx$.
These densities are nearly indistinguishable and, in fact,
their total variation distance is $e^{-1/\varepsilon}$.}
\label{fig::picture3}
\end{figure}

Figure \ref{fig::picture3}
shows the type of least favorable pair
we used for manifold learning.
The top two plots do not show the densities; rather they show the
support of the densities.
The distribution $g_0$ is uniform on the circle in the top left plot.
The distribution $g_1$ is uniform on the perturbed circle in the top
right plot.
The Hausdorff distance between the supports of densities is $\varepsilon$.
The bottom plots are
$q_0 = \int\phi(y-x) g_0(x)\, dx$ and
$q_1 = \int\phi(y-x) g_1(x) \,dx$.
Again, these\vadjust{\goodbreak} densities are nearly indistinguishable, and, in fact,
their total variation distance is $e^{-1/\varepsilon}$.
In this case, however, $g_0$ and $g_1$
have different supports.

\subsection{Relationship to regression with measurement error}
We can also relate the manifold estimation problem
with nonparametric regression with measurement error.
Suppose that
%
\begin{eqnarray}\label{eq::measerr}
U_i &=& X_i + Z_{2i}, \\
\nonumber
Y_i &=& m(X_i) + Z_{1i},\vadjust{\goodbreak}
\end{eqnarray}
and we want to estimate the regression function $m$.
If we observe
$(X_1,Y_1), \ldots,\allowbreak (X_n,Y_n)$, then this is a~standard
nonparametric regression problem.
But if we only observe
$(U_1,Y_1),\ldots, (U_n,Y_n)$, then this is
the usual nonparametric regression with measurement error problem.
The rates of convergence are similar to deconvolution.
Indeed, \citet{FanTruong} have an argument that
converts
nonparametric regression with measurement error
into a~density deconvolution problem.
Let us see how this related to manifold learning.

Suppose that
$D=2$ and $d=1$.
Futher, suppose that the manifold is
\textit{function-like},
meaning that the manifold is a~curve of the form
$M=\{ (u,\break m(u))\dvtx u\in\mathbb{R}\}$ for some function $m$.
Then each $Y_i$ can be written in the form
\[
Y_i =
\pmatrix{ Y_{i1} \cr Y_{i2}}
  =
\pmatrix{ U_i \cr m(U_i)
}+
\pmatrix{ Z_{1i} \cr Z_{2i}
}
\]
which is exactly of the form
(\ref{eq::measerr}).
Let $\mathcal{Q}$ be all such distributions obtained this way with
$|m''(u)| \leq1/\kappa$.
However, this only holds when the manifold
has the function-like form.
Moreover, the lower bound argument in \citet{FanTruong}
cannot directly be transferred to the manifold setting as we now explain.

In our lower bound proof, we defined a~least favorable
pair $q_0$ and $q_1$
for the distribution of $Y$ as follows.
Take
$M_0=\{ (u,0)\dvtx u\in\mathbb{R}\}$ and
$M_1=\{ (u,m(u))\dvtx u\in\mathbb{R}\}$.
[In fact, we used $(u,m(u))$ and
$(u,-m(u))$,
but the present discussion is clearer if we use
$(u,0)$ and $(u,m(u))$.]
Let $Y = (Y_1,Y_2)$.
For $M_0$, the distribution $q_0$ for $Y$ is based on
\[
 \pmatrix{ Y_1 \cr Y_2
} =
\pmatrix{ U \cr 0
} +
\pmatrix{ Z_1 \cr Z_2
}.
\]
The density of $(U,Y_2)$ is
$f_0(u,y_2)= \zeta(u)\phi(y_2)$
where $\zeta$ is some density for~$U$.
Then
\[
q_0(y_1,y_2) = f_0 \star\Phi=
\int f_0(y_1-Z_1,y_2) \,d\Phi(z_1),
\]
where the convolution symbol here and in what follows, refers to
convolution only over $U + Z_1$.

Now let
$q_1(y_1,y_2)$ denote the distribution of $Y$ in the model
\[
 \pmatrix{ Y_1 \cr Y_2
} =
 \pmatrix{ U \cr m(U)
} +
 \pmatrix{ Z_1 \cr Z_2
}.
\]
This generates the least favorable pair $q_0$, $q_1$ used in our proof
(restricted to this special case).

The least favorable pair used by Fan and Truong is different in a~subtle way.
The first distribution $q_0$ is the same.
The second, which we will denote~$w_1$, is constructed as follows.
Let
\[
w_1(y_1,y_2) = f_1 \star\Phi,
\]
where
the convolution is only over $U$,
\[
f_1(\xi,y_2)= f_0(\xi,y_2) + \gamma
H\bigl(\xi/\sqrt{\gamma}\bigr)h_0(y_2),
\]
where
$f_1(\xi)=g(\xi)$,
$\gamma H(\xi/\sqrt{\gamma})/g(\xi) = b(\xi)$,
$H$ is a~perturbation function such as a~cosine,
and
$h_0$ is chosen so that
$\int h_0(y_2) \,dy_2 =0$ and
$\int y_2 h_0(y_2) \,dy_2 =1$.
Now we show that $w_1(y_1,y_2) \neq q_1(y_1,y_2)$.
In fact, $w_1$ is not in $\mathcal{Q}$.
Note that
\[
w_1(y_1,y_2) = f_1 \star\Phi=
q_0(y_1,y_2) + \gamma h_0(y_2)
\int H \biggl(\frac{y_1-z_1}{\sqrt{\gamma}} \biggr)\,d\Phi(z_1).
\]
Now,
\[
q_1(y_2|u) = \phi\bigl(y_2 - m(u)\bigr),
\]
but
\[
f_1(y_2|u) = \frac{f_1(y_2,u)}{f_1(u)} = \phi(y_2) + m(u)h_0(y_2).
\]
These both have mean $m(u)$ but the distributions are different.
Indeed, the marginals
$w_1(y_2)$ and $q_1(y_2)$ are different.
In fact,
\[
w_1(y_2) = q_0(y_2) + c h_0(y_2)
\]
for some $c$.
This is not in our class because it is not of the form
$\phi(y_2-m(u))$.
Hence,
$w_1$ is not in our class $\mathcal{Q}$:
it does not correspond to drawing a~point on a~manifold and adding noise.

The point is that
manifold learning reduces to nonparametric regression
with errors only in the special case that the manifold is
function-like.
And even in this case, the proofs of the bounds
are somewhat different than the usual proofs.

\section{Discussion}

The purpose of this paper is to establish
minimax bounds on estimating
manifolds.
The estimators used to prove the upper bounds are
theoretical constructions for the purposes of the proofs.
They are not practical estimators.

There is a~large literature on methodology for estimating manifolds.
However, these estimators are not likely to be optimal except under
stringent conditions.
In current work we are trying to bridge the gap between the theory
and the methodology.

Probably the most realistic noise condition is the additive model.
In this case, we are dealing with a~singular deconvolution problem.
The upper bound used deconvolution techniques.
Such methods require that the noise distribution is known
(or is at least restricted to some narrow class of distributions).
This seems unrealistic in real problems.
A~more realistic goal is to estimate some
proxy manifold $M^*$ that, in some sense, approximates $M$.
We are currently working on such techniques.


%

\printaddresses


\begin{thebibliography}{20}

\bibitem[\protect\citeauthoryear{Arias-Castro et~al.}{2005}]{dots}
%
\begin{barticle}[mr]
\bauthor{\bsnm{Arias-Castro},~\bfnm{Ery}\binits{E.}},
\bauthor{\bsnm{Donoho},~\bfnm{David~L.}\binits{D.~L.}},
\bauthor{\bsnm{Huo},~\bfnm{Xiaoming}\binits{X.}} \AND
\bauthor{\bsnm{Tovey},~\bfnm{Craig~A.}\binits{C.~A.}}
(\byear{2005}).
\btitle{Connect the dots: How many random points can a~regular curve pass
through?}
\bjournal{Adv. in Appl. Probab.}
\bvolume{37}
\bpages{571--603}.
\bid{doi={10.1239/aap/1127483737}, issn={0001-8678}, mr={2156550}}
\bptok{imsref}%
\end{barticle}
%
\endbibitem

\bibitem[\protect\citeauthoryear{Baraniuk and Wakin}{2009}]{baraniuk}
%
\begin{barticle}[mr]
\bauthor{\bsnm{Baraniuk},~\bfnm{Richard~G.}\binits{R.~G.}} \AND
\bauthor{\bsnm{Wakin},~\bfnm{Michael~B.}\binits{M.~B.}}
(\byear{2009}).
\btitle{Random projections of smooth manifolds}.
\bjournal{Found. Comput. Math.}
\bvolume{9}
\bpages{51--77}.
\bid{doi={10.1007/s10208-007-9011-z}, issn={1615-3375}, mr={2472287}}
\bptnote{check year}%
\bptok{imsref}%
\end{barticle}
%
\endbibitem

\bibitem[\protect\citeauthoryear{Bousquet, Boucheron and Lugosi}{2004}]{BBL}
%
\begin{bmisc}[auto:STB|2012/04/12|05:18:16]
\bauthor{\bsnm{Bousquet},~\bfnm{O.}\binits{O.}},
\bauthor{\bsnm{Boucheron},~\bfnm{S.}\binits{S.}} \AND
\bauthor{\bsnm{Lugosi},~\bfnm{G.}\binits{G.}}
(\byear{2004}).
\bhowpublished{Introduction to statistical learning theory. \textit{Machine
Learning} \textbf{3176} 169--207}.
\bptok{imsref}%
\end{bmisc}
%
\endbibitem

\bibitem[\protect\citeauthoryear{Caillerie et~al.}{2011}]{CCDM}
%
\begin{barticle}[mr]
\bauthor{\bsnm{Caillerie},~\bfnm{Claire}\binits{C.}},
\bauthor{\bsnm{Chazal},~\bfnm{Fr{\'e}d{\'e}ric}\binits{F.}},
\bauthor{\bsnm{Dedecker},~\bfnm{J{\'e}r{\^o}me}\binits{J.}} \AND
\bauthor{\bsnm{Michel},~\bfnm{Bertrand}\binits{B.}}
(\byear{2011}).
\btitle{Deconvolution for the {W}asserstein metric and geometric inference}.
\bjournal{Electron. J. Stat.}
\bvolume{5}
\bpages{1394--1423}.
\bid{doi={10.1214/11-EJS646}, issn={1935-7524}, mr={2851684}}
\bptok{imsref}%
\end{barticle}
%
\endbibitem

\bibitem[\protect\citeauthoryear{Chaudhuri and Dasgupta}{2010}]{SanjoyTree}
%
\begin{bmisc}[auto:STB|2012/04/12|05:18:16]
\bauthor{\bsnm{Chaudhuri},~\bfnm{K.}\binits{K.}} \AND
\bauthor{\bsnm{Dasgupta},~\bfnm{S.}\binits{S.}}
(\byear{2010}).
\bhowpublished{Rates of convergence for the cluster tree. In \textit{Advances in Neural Information Processing Systems 23}
(J. Lafferty, C. K. I. Williams, J. Shawe-Taylor, R.~S. Zemel and A. Culotta) 343--351}.
\bptok{imsref}%
\end{bmisc}
%
\endbibitem

\bibitem[\protect\citeauthoryear{Dey}{2007}]{Dey}
%
\begin{bbook}[mr]
\bauthor{\bsnm{Dey},~\bfnm{Tamal~K.}\binits{T.~K.}}
(\byear{2007}).
\btitle{Curve and Surface Reconstruction: Algorithms with Mathematical
Analysis}.
\bseries{Cambridge Monographs on Applied and Computational Mathematics}
\bvolume{23}.
\bpublisher{Cambridge Univ. Press}, \baddress{Cambridge}.
\bid{mr={2267420}}
\bptnote{check year}%
\bptok{imsref}%
\end{bbook}
%
\endbibitem

\bibitem[\protect\citeauthoryear{Fan}{1991}]{Fan}
%
\begin{barticle}[mr]
\bauthor{\bsnm{Fan},~\bfnm{Jianqing}\binits{J.}}
(\byear{1991}).
\btitle{On the optimal rates of convergence for nonparametric deconvolution
problems}.
\bjournal{Ann. Statist.}
\bvolume{19}
\bpages{1257--1272}.
\bid{doi={10.1214/aos/1176348248}, issn={0090-5364}, mr={1126324}}
\bptok{imsref}%
\end{barticle}
%
\endbibitem

\bibitem[\protect\citeauthoryear{Fan and Truong}{1993}]{FanTruong}
%
\begin{barticle}[mr]
\bauthor{\bsnm{Fan},~\bfnm{Jianqing}\binits{J.}} \AND
\bauthor{\bsnm{Truong},~\bfnm{Young~K.}\binits{Y.~K.}}
(\byear{1993}).
\btitle{Nonparametric regression with errors in variables}.
\bjournal{Ann. Statist.}
\bvolume{21}
\bpages{1900--1925}.
\bid{doi={10.1214/aos/1176349402}, issn={0090-5364}, mr={1245773}}
\bptok{imsref}%
\end{barticle}
%
\endbibitem

\bibitem[\protect\citeauthoryear{Federer}{1959}]{federer}
%
\begin{barticle}[mr]
\bauthor{\bsnm{Federer},~\bfnm{Herbert}\binits{H.}}
(\byear{1959}).
\btitle{Curvature measures}.
\bjournal{Trans. Amer. Math. Soc.}
\bvolume{93}
\bpages{418--491}.
\bid{issn={0002-9947}, mr={0110078}}
\bptok{imsref}%
\end{barticle}
%
\endbibitem

\bibitem[\protect\citeauthoryear{Genovese et~al.}{2010}]{us::2010}
%
\begin{bmisc}[auto:STB|2012/04/12|05:18:16]
\bauthor{\bsnm{Genovese},~\bfnm{C.~R.}\binits{C.~R.}},
\bauthor{\bsnm{Perone-Pacifico},~\bfnm{M.}\binits{M.}},
\bauthor{\bsnm{Verdinelli},~\bfnm{I.}\binits{I.}} \AND
\bauthor{\bsnm{Wasserman},~\bfnm{L.}\binits{L.}}
(\byear{2010}).
\bhowpublished{Minimax manifold estimation. Available at
arXiv:\arxivurl{1007.0549}}.
\bptok{imsref}%
\end{bmisc}
%
\endbibitem


\bibitem[\protect\citeauthoryear{Koltchinskii}{2000}]{Koltchinskii}
%
\begin{barticle}[mr]
\bauthor{\bsnm{Koltchinskii},~\bfnm{V.~I.}\binits{V.~I.}}
(\byear{2000}).
\btitle{Empirical geometry of multivariate data: A~deconvolution approach}.
\bjournal{Ann. Statist.}
\bvolume{28}
\bpages{591--629}.
\bid{doi={10.1214/aos/1016218232}, issn={0090-5364}, mr={1790011}}
\bptok{imsref}%
\end{barticle}
%
\endbibitem

\bibitem[\protect\citeauthoryear{Meister}{2006}]{Meister20061702}
%
\begin{barticle}[mr]
\bauthor{\bsnm{Meister},~\bfnm{Alexander}\binits{A.}}
(\byear{2006}).
\btitle{Estimating the support of multivariate densities under measurement
error}.
\bjournal{J.~Multivariate Anal.}
\bvolume{97}
\bpages{1702--1717}.
\bid{doi={10.1016/j.jmva.2005.04.004}, issn={0047-259X}, mr={2298884}}
\bptok{imsref}%
\end{barticle}
%
\endbibitem

\bibitem[\protect\citeauthoryear{Niyogi, Smale and Weinberger}{2008}]{smale}
%
\begin{barticle}[mr]
\bauthor{\bsnm{Niyogi},~\bfnm{Partha}\binits{P.}},
\bauthor{\bsnm{Smale},~\bfnm{Stephen}\binits{S.}} \AND
\bauthor{\bsnm{Weinberger},~\bfnm{Shmuel}\binits{S.}}
(\byear{2008}).
\btitle{Finding the homology of submanifolds with high confidence from random
samples}.
\bjournal{Discrete Comput. Geom.}
\bvolume{39}
\bpages{419--441}.
\bid{doi={10.1007/s00454-008-9053-2}, issn={0179-5376}, mr={2383768}}
\bptnote{check year}%
\bptok{imsref}%
\end{barticle}
%
\endbibitem

\bibitem[\protect\citeauthoryear{Niyogi, Smale and
Weinberger}{2011}]{smale2}
%
\begin{bmisc}[auto:STB|2012/04/12|05:18:16]
\bauthor{\bsnm{Niyogi},~\bfnm{P.}\binits{P.}},
\bauthor{\bsnm{Smale},~\bfnm{S.}\binits{S.}} \AND
\bauthor{\bsnm{Weinberger},~\bfnm{S.}\binits{S.}}
(\byear{2011}).
\bhowpublished{A~topological view of unsupervised learning from noisy data.
\textit{SIAM J. Comput.} \textbf{40} 646--663}.
\bptok{imsref}%
\end{bmisc}
%
\endbibitem

\bibitem[\protect\citeauthoryear{Ozertem and Erdogmus}{2011}]{Principal}
%
\begin{barticle}[mr]
\bauthor{\bsnm{Ozertem},~\bfnm{Umut}\binits{U.}} \AND
\bauthor{\bsnm{Erdogmus},~\bfnm{Deniz}\binits{D.}}
(\byear{2011}).
\btitle{Locally defined principal curves and surfaces}.
\bjournal{J. Mach. Learn. Res.}
\bvolume{12}
\bpages{1249--1286}.
\bid{issn={1532-4435}, mr={2804600}}
\bptok{imsref}%
\end{barticle}
%
\endbibitem

\bibitem[\protect\citeauthoryear{Stefanski}{1990}]{Stefanski1990229}
%
\begin{barticle}[mr]
\bauthor{\bsnm{Stefanski},~\bfnm{Leonard~A.}\binits{L.~A.}}
(\byear{1990}).
\btitle{Rates of convergence of some estimators in a~class of deconvolution
problems}.
\bjournal{Statist. Probab. Lett.}
\bvolume{9}
\bpages{229--235}.
\bid{doi={10.1016/0167-7152(90)90061-B}, issn={0167-7152}, mr={1045189}}
\bptok{imsref}%
\end{barticle}
%
\endbibitem

\bibitem[\protect\citeauthoryear{Yu}{1997}]{binyu}
%
\begin{bincollection}[mr]
\bauthor{\bsnm{Yu},~\bfnm{Bin}\binits{B.}}
(\byear{1997}).
\btitle{Assouad, {F}ano, and {L}e {C}am}.
In \bbooktitle{Festschrift for {L}ucien {L}e {C}am}
\bpages{423--435}.
\bpublisher{Springer}, \baddress{New York}.
\bid{mr={1462963}}
\bptok{imsref}%
\end{bincollection}
%
\endbibitem

\bibitem[\protect\citeauthoryear{Yukich}{1985}]{Yukich1985245}
%
\begin{barticle}[mr]
\bauthor{\bsnm{Yukich},~\bfnm{J.~E.}\binits{J.~E.}}
(\byear{1985}).
\btitle{Laws of large numbers for classes of functions}.
\bjournal{J. Multivariate Anal.}
\bvolume{17}
\bpages{245--260}.
\bid{doi={10.1016/0047-259X(85)90083-1}, issn={0047-259X}, mr={0813235}}
\bptok{imsref}%
\end{barticle}
%
\endbibitem

\end{thebibliography}
\end{document}